\documentclass[11pt]{article}
 \usepackage[psamsfonts]{amssymb}
 \usepackage[applemac]{inputenc}

\newcommand{\bk}{l\!k}

\title{An asymptotic formula for the ranks of    the homotopy groups of a finite complex}

\author{Yves Felix,   Steve Halperin   and Jean-Claude Thomas}

\begin{document}
\maketitle

    \begin{abstract}
      Let $X$ be a finite simply connected
CW complex of dimension $n$. The loop space homology $H_*(\Omega
X;\mathbb Q)$ is the universal enveloping algebra of a graded Lie
algebra $L_X$ isomorphic with $\pi_{*-1} (X)\otimes \mathbb Q$.
Let $Q_X \subset L_X$ be a minimal generating subspace, and set
$\alpha = \limsup_i \frac{\log\mbox{\scriptsize rk}\,
\pi_i(X)}{i}$.  Theorem:\\ If $\mbox{dim}\, L_X = \infty$ and
$\limsup (\mbox{dim}\,( Q_X)_k)^{1/k} < \limsup (\mbox{dim}\,
(L_X)_k)^{1/k}\,,$  then  $$\sum_{i=1}^{n-1}
\mbox{rk}\,\pi_{k+i}(X) = e^{(\alpha + \varepsilon_k)k}\,,
\hspace{1cm} \mbox{where}\, \varepsilon_k \to 0\,\,\mbox{as}\,
k\to \infty\,.$$
   In
particular $\displaystyle\sum_{i=1}^{n-1} \, \mbox{rk}\,
\pi_{k+i}(X)$ grows exponentially in $k$.
  \end{abstract}

    \vspace{5mm}\noindent {\bf AMS Classification} : 55P35, 55P62,
    17B70

    \vspace{2mm}\noindent {\bf Key words} : Homotopy Lie algebra,
    graded Lie algebra, exponential growth

\section{Introduction}

Suppose $X$ is a finite simply connected CW complex of dimension
$n$. The homotopy groups of $X$ then have the form
$$\pi_i(X) = \mathbb Z^{\rho_i} \oplus T_i\,, $$
where $T_i$ is a finite abelian group   and
  $\rho_i = $ rk$\, \pi_i(X)$ is finite. It is known  \cite{FHT}
   that either $\pi_i(X) = T_i$, $i\geq 2n$ ($X$ is {\it rationally elliptic})
    or else for all $k\geq 1$, $\sum_{i=1}^{n-1} $ rk$\, \pi_{k+i}(X) >
    0$.
In this case $X$ is called {\it rationally hyperbolic}.

In [7] it is shown that in the rationally hyperbolic case $\sum_{i
  = 1}^{n-1} $ rk$\, \pi_{k+i}(X) $ grows faster than any polynomial in
$k$. Here we show that with an additional hypothesis this sum
grows exponentially in $k$ and, in fact, setting $\alpha = \limsup
\frac{\log \mbox{\scriptsize rk}\, \pi_i}{i}$ we have $$ \sum_{i 1}^{n-1} \mbox{rk}\, \pi_{k+i} =e^{(\alpha +\varepsilon_k)k}\,,
\hspace{2cm}\mbox{where}\, \varepsilon_k\to 0\,\, \mbox{as }\,
k\to \infty\,.$$

In   subsequent papers we will identify a large class of spaces
for which the additional hypothesis holds: in fact it may well
hold for all finite simply connected CW complexes.

 Note that
rk$\,\pi_i(X) = \mbox{dim}\, \pi_i(X)\otimes \mathbb Q$. Thus we
work more generally with simply connected spaces $X$ such that
each $H_i(X;\mathbb Q)$ is finite dimensional. In this case $
\mbox{dim}\, \pi_i(X)\otimes \mathbb Q$ is also finite for each
$i$. On the other hand, a theorem of Milnor-Moore-Cartan-Serre
asserts that the loop space homology $H_*(\Omega X;\mathbb Q)$ is
the universal enveloping algebra of a graded Lie algebra $L_X$ and
that the Hurewicz homomorphism is an isomorphism $\pi_*(\Omega
X)\otimes \mathbb Q \stackrel{\cong}{\longrightarrow} L_X$. Since
there are natural isomorphisms $\pi_i(X) \cong \pi_{i-1}(\Omega
X)$ it follows that
  the results above
can be phrased in terms of the integers $\mbox{dim}\,( L_X)_i$.

For any graded vector space $V$ concentrated in positive degrees
we define the {\it  logarithmic index } of $V$ by $$\mbox{log
index}\, V = \limsup \frac{\log \mbox{dim}\, V_k}{k}\,.$$ In
(\cite{IHES}, \cite{FHT}) it is shown that if $X$ (simply
connected) has finite Lusternik-Schnirelmann category (in
particular, if $X$ is a finite CW complex) and if $X$ is
rationally hyperbolic then $$\mbox{log index}\, L_X > 0 \,.$$

Now let $Q_X$ denote a minimal generating subspace for the Lie
algebra $L_X$, and let $\alpha$ denote log index$\, L_X$.

\vspace{3mm}\noindent {\bf Theorem 1.}  {\sl Let $X$ be a simply
connected topological space with finite dimensional rational
homology concentrated in degrees $\leq n$. Suppose  $\mbox{log
index}\, Q_X <$ log index$\, L_X$, then $$\sum_{i=1}^{n-1}
\mbox{dim}\, (L_X)_{k+i} = e^{(\alpha +\varepsilon_k)k}\,,
\hspace{1cm} \mbox{where}\, \varepsilon_k\to 0\,\, \mbox{as} \,
k\to\infty\,.$$ In particular this sum grows exponentially in $k$.
}

\vspace{3mm}\noindent {\bf Theorem 2.}   {\sl Let $X$ be a simply
connected topological space with finite dimensional rational
homology in each degree, and finite   Lusternik-Schnirelmann
category. Suppose log index$\, Q_X<$ log index$\, L_X<\infty$.
Then for some   $d>0$, $$\sum_{i=1}^{d-1} \mbox{dim}\, (L_X)_{k+i}
= e^{(\alpha +\varepsilon_k)k}\,, \hspace{1cm} \mbox{where}\,
\varepsilon_k\to 0\,\, \mbox{as} \, k\to\infty\,.$$   In
particular, $\sum_{i=1}^{d-1}\mbox{dim}\, (L_X)_{k+i}$ grows
exponentially in $k$. }

\vspace{3mm}\noindent {\bf Corollary:}  The conclusion of Theorem
1 holds for finite simply connected CW complexes $X$ for which
$L_X$ is infinite, but finitely generated. The conclusion of
Theorem 2 holds for simply connected spaces of finite LS category
and finite rational Betti numbers provided that $L_X$ is infinite,
but  finitely generated, and, log index$\, L_X
< \infty$.

\vspace{3mm} Recall that the depth of a graded  Lie
algebra $L$ is the least $m$ (or $\infty$) such that
$\mbox{Ext}_{UL}^m(\mathbb Q, UL)\neq 0$.

The key  ingredients in the proofs of Theorems 1 and 2 are
\begin{enumerate}
\item[$\bullet$] A growth condition for $L_X$ established in
(\cite{IHES},\cite{FHT})\item[$\bullet$] The fact that depth $L_X
< \infty$, established in (\cite{rad},\cite{FHT})
\end{enumerate}

We shall use Lie algebra arguments in Theorem 3 below to deduce
the conclusion of Theorems 1 and 2 from these ingredients, and
then deduce Theorems 1 and 2.

Theorems 1 and 2 may be compared with the results in \cite{Roos}
and in \cite{ENS} that assert  for $n$-dimensional finite CW
complexes (respectively for simply connected spaces with finite
type rational homology and   finite Lusternik-Schnirelmann
category) that $ \sum_{i=1}^{n-1}\mbox{dim}\, (L_X)_{k+i}$ (resp.
$\sum_{i=1}^{d-1}\mbox{dim}\, (L_X)_{k+i}$) grows faster than any
polynomial in $k$. These results use only the fact that $L_X$ has
finite depth, and require no hypothesis on $Q_X$.

The hypothesis on $Q_X$ in Theorems 1 and 2 may be restated as
requiring that the formal series \hspace{3mm} $\sum_q\mbox{dim}\,
\mbox{Tor}_{1,q}^{UL_X}\, z^q$ \hspace{3mm} have a radius of
convergence strictly greater than that of the formal series
$\sum_q \mbox{dim}\, (UL_X)_q \,z^q$. Lambrechts \cite{La1} has
proved a much stronger result under the hypothesis that the formal
series $\sum_q\left(\sum_p (-1)^p \mbox{dim}\,
\mbox{Tor}_{p,q}\right) z^q$ has a radius of convergence strictly
larger than that of $\sum_q \mbox{dim}\, (UL_X)_q\, z^q$.

\section{Lie algebras}

In this section we work over any ground field $\bk$ of
characteristic different from $2$;  graded Lie algebras $L$ are
defined as in {\cite{Roos} and, in particular, are assumed to
satisfy $[x,[x,x]] = 0$, $x\in L_{\mbox{\scriptsize  odd}}$ (This
follows from the Jacobi identity except  when char $\bk = 3$).

A graded Lie algebra  $L$ is connected and of finite type if
\begin{center} $L = \{ L_i\}_{i\geq 1}$ and each $L_i$ is finite
dimensional. \end{center} We shall refer to these as {\it  cft Lie
algebras}.  The minimal generating subspaces $Q$ of a cft Lie
algebra $L$  are those   subspaces $Q$ for which $Q \to L/[L,L]$
is a linear isomorphism.

A {\it  growth sequence} for a cft Lie algebra $L$ is a sequence
$(r_i)$ such that $r_i \to \infty$ and
$$\lim_{i\to \infty} \frac{\log \mbox{dim}\, L_{r_i}}{r_i} = \mbox{log
  index}\, L\,.$$
A {\it  quasi-geometric sequence}  $(\ell_i)$ is a sequence such that
for some integer $m$, $\ell_i< \ell_{i+1}\leq m\ell_i$, all $i$; if
additionally $(\ell_i)$ is a growth sequence then it is a {\it
  quasi-geometric growth sequence}.

Of particular interest here are the growth conditions $$0 <
\mbox{log index}\, L < \infty\,; \eqno{(A.1)}$$ and $$\mbox{log
index}\, Q < \mbox{log index} L\,. \eqno{(A.2)}$$

\vspace{3mm}\noindent{\bf Proposition.}  {\sl Let $L$ be a cft Lie
algebra satisfying   (A.1) and (A.2), and assume
$L$ has a quasi-geometric growth sequence $(r_j)$. Then any sequence
 $(s_i)$  such that
$s_i \to \infty$ has a subsequence $(s_{i_j})$ for which there are
growth sequences  $(t_j)$ and $(p_j)$ such that $$t_j\leq s_{i_j} <
p_j\hspace{1cm}\mbox{and} \hspace{4mm} p_j/t_j \to 1\,.$$}

\vspace{4mm}\noindent {\sl Proof.}  First note that because of (A.2),
$$\mbox{log index}\, L/Q = \mbox{log index}\, L$$
Now  adopt the following
notation, for $i\geq 1$ :  $$ \left. \begin{array}{l}
\mbox{log index}\, L/Q = \alpha\\ \mbox{dim}\, L_i e^{(\alpha + \varepsilon_i)i}\\ \mbox{dim}\, Q_i = e^{(\alpha +
\sigma_i)i}\\ \mbox{dim}\, L_i/Q_i = e^{(\alpha +\tau_i)i}
\\\mbox{dim}\, (UL)_i = e^{(\alpha + \delta_i)i}
\end{array}
\right\} \eqno{(1)}$$ Then because of (A.1) and  (A.2), $0<\alpha <
\infty$, and  $\limsup
\varepsilon_i = 0$. Moreover, by a result of Babenko (\cite{Ba},
\cite{FHT}),
$\mbox{log index}\, UL = \mbox{log index}\, L$,
 and so $\limsup (\delta_i) = 0$. Finally  (A.2)
implies that $\limsup (\sigma_i) = \sigma < 0$, and that $\tau_{q_i}
\to 0$ as $i\to \infty$.

Next, since $r_j$ is a quasi-geometric sequence, for some fixed $m$ we
have $r_j<r_{j+1} \leq mr_j$, all $j$. It follows that for each $s_i$
in our sequence we may choose $q_i$ in the sequence $(r_j)$ so that
 $$s_i<q_i\leq ms_i\,.\eqno{(2)}$$ Thus $q_i \to
\infty$ as $i\to \infty$. The adjoint representation of $UL$ in
$L$ defines surjections $$ \displaystyle\bigoplus_{(\ell ,k,t)\in
{\cal J}_i} (UL)_{\ell}\otimes Q_k\otimes L_t \twoheadrightarrow
\displaystyle\bigoplus_{({\ell},k,t)\in{\cal J}_i}
(UL)_{\ell}\otimes [Q_k,L_t] \twoheadrightarrow
L_{q_i}/Q_{q_i}\,,\eqno{(3)}$$ where ${\cal J}_i$ consists of
those triples for which ${\ell}+k+t = q_i$ and $t\leq s_i < t+k$.

\noindent Thus $$3(q_i+1)\renewcommand{\arraystretch}{0.8}
\begin{array}[t]{l}\mbox{max}\\
{\scriptstyle (\ell ,k,t)\in {\cal J}_i} \end{array}
\renewcommand{\arraystretch}{1} \mbox{dim}\, (UL)_{\ell} \,\,
\mbox{dim}\, [Q_k,L_t] \geq \mbox{dim}\, L_{q_i}/Q_{q_i}\,.$$ For
each $s_i$ we may therefore choose $(\ell_i,k_i,t_i)\in {\cal
J}_i$ so that $$3(q_i+1) \, \mbox{dim}\, (UL)_{\ell_i}\,
\mbox{dim}\, [Q_{k_i},L_{t_i}] \geq \mbox{dim}\,
L_{q_i}/Q_{q_i}\,.\eqno{(4)}$$

\vspace{3mm}\noindent {\bf Lemma 1.}  {\sl $(k_i + t_i) \geq
\frac{1}{m} \,q_i$, and  $(k_i+t_i) \geq \frac{1}{m-1}\ell_i$. In
particular, $k_i+t_i \to \infty$ as $i\to \infty$.}

\vspace{3mm}\noindent {\sl Proof.}  It follows from (2) that $k_i
+ t_i \geq s_i \geq \frac{1}{m}\, q_i$. Since $q_i = k_i +t_i+\ell_i$,
$(k_i+t_i)\geq \frac{1}{m-1}\ell_i$. Since $s_i \to \infty$ so
does $k_i + t_i$. \hfill $\square$

\vspace{7mm}

Next, define $\varepsilon (k_i,t_i)$ by $$\mbox{dim}\,
[Q_{k_i},L_{t_i}] = e^{[\alpha + \varepsilon
(k_i,t_i)](k_i+t_i)}\,.$$

\vspace{3mm}\noindent {\bf Lemma 2.} {\sl $\varepsilon (k_i,t_i)
\to 0$ and $\varepsilon_{k_i+t_i} \to 0$ as $i\to \infty$.}

\vspace{3mm}\noindent {\bf Proof.}  Define $\epsilon (q_i)$ by
$3(q_i + 1) = e^{\varepsilon (q_i)q_i}$. Then (4) reduces to
$$\varepsilon (q_i)q_i + (\alpha + \delta_{\ell_i}) \ell_i + [\alpha +
\varepsilon (k_i,t_i)](k_i + t_i) \geq (\alpha + \tau_{q_i}
) q_i\,.$$ Since $\ell_i + k_i + t_i = q_i$, $$\varepsilon (q_i)
\frac{q_i}{k_i+t_i} + \delta_{\ell_i} \frac{\ell_i}{k_i + t_i} +
\varepsilon (k_i,t_i) \geq \tau_{q_i}
\frac{q_i}{k_i+t_i}\,.$$ Now as $i\to \infty$, $q_i \geq s_i \to
\infty$. Thus $\varepsilon (q_i) \to 0$. Moreover, since $q_i$
belongs to a growth sequence, $\varepsilon_{q_i} \to 0$ and
$\tau_{q_i} \to 0$.  Use Lemma
1 to conclude that $$\varepsilon (q_i) \frac{q_i}{k_i+t_i} \to 0
\hspace{1cm} \mbox{and} \hspace{5mm} \tau_{q_i}
\frac{q_i}{k_i+t_i}\to 0\,,$$ and hence $$\liminf \varepsilon
(k_i,t_i) \geq - \limsup \delta_{\ell_i} \frac{\ell_i}{k_i+t_i}\,.$$

Next, since $\limsup\, \delta_{\ell_i} = 0$ and $\frac{\ell_i}{k_i+t_i}
\leq m-1$, it follows that
$$- \limsup \delta_{\ell_i} \frac{\ell_i}{k_i+t_i} = 0\,.$$

Finally, $[Q_{k_i},L_{t_i}]$ embeds in $L_{k_i+t_i} / Q_{k_i +
t_i}$, and it follows that $$\varepsilon_{k_i+t_i} \geq
\varepsilon (k_i, t_i)\,.$$ But $k_i+t_i \to \infty$ and so
 $\limsup\, \varepsilon_{k_i+t_i} \leq 0$. This, together
with $\liminf \, \varepsilon (k_i,t_i) \geq 0$,
completes the proof of the lemma. \hfill $\square$

\vspace{1cm} Next, since $Q_{k_i} \otimes L_{t_i} \to [Q_{k_i},
L_{t_i}]$ is surjective we have $$\sigma_{k_i}k_i +
\varepsilon_{t_i}t_i \geq \varepsilon (k_i,t_i)
(k_i+t_i)\eqno{(5)}$$

\vspace{3mm}\noindent {\bf Lemma 3.}  {\sl $t_i/k_i \to \infty$ and $t_i \to \infty$ as $i
\to \infty$.}

\vspace{3mm}\noindent {\sl Proof.} If $t_i/k_i$ does not converge to
$\infty$ then  we would have
$t_{i_\nu}/k_{i_\nu} \leq T$ for some subsequence $(s_{i_\nu})$.  But
$$\sigma_{k_{i_\nu}} \geq -\varepsilon_{t_{i_\nu}}
\frac{t_{i_\nu}}{k_{i_\nu}} + \varepsilon (k_{i_\nu}, t_{i_\nu})
(1 + t_{i_\nu}/k_{i_\nu})\,.$$ Since (Lemma 1) $k_{i_\nu} + t_{i_\nu} \to
\infty$ and (Lemma 2) $\varepsilon (k_{i_\nu}, t_{i_\nu}) \to 0$,
the $\liminf$ of the  right hand side of this equation would be $\geq 0$. Hence $\limsup \,
\sigma_{k_{i_\nu}} \geq 0$, which would contradict $\limsup \,
\sigma_i < 0$. Finally, since $t_i+k_i \to \infty$ it follows that
$t_i\to \infty$. \hfill $\square$

\vspace{3mm}\noindent{\bf Lemma 4.}  {\sl Write $k_i \lambda_it_i$. Then, $$\lambda_i \to 0\hspace{4mm}\mbox{and }
\,\, \varepsilon_{t_i} \to 0\hspace{15mm}\mbox{as}\hspace{3mm}
i\to \infty\,.$$ }

\vspace{3mm}\noindent {\sl Proof.} Lemma 3 asserts that $\lambda_i \to
0$. Rewrite equation (5) as
$$\varepsilon_{t_i} \geq (\varepsilon (k_i,t_i)  -
\sigma_{k_i})\lambda_i + \varepsilon (k_i,t_i)\,. \eqno{(6)}$$
Since $\limsup \sigma_j$ is finite, the $\sigma_j$ are bounded
above. Thus $-\sigma_j \geq A$, some constant $A$. Since $\lambda_i
\to 0$,  $\liminf
(-\sigma_{k_i}\lambda_i)\geq 0$. Since (Lemma 2) $ \varepsilon (k_i,t_i)
\to 0$ it follows from (6) that $\liminf \varepsilon_{t_i} = 0$. But
$\limsup \varepsilon_{t_i} \leq \limsup \varepsilon_i = 0$ and so
$\varepsilon_{t_i} \to 0$.
 \hfill
$\square$.

\vspace{4mm} The lemmas above establish the Proposition. Simply
set $p_i = t_i + k_i$ and note that $t_i \to \infty $ (Lemma 3),
$p_i \to \infty$ (Lemma 1), $p_i/t_i = 1 + \lambda_i \to 1$ (Lemma
4). Furthermore $\varepsilon_{t_i} \to 0$ (Lemma 4) and
$\varepsilon_{p_i} \to 0$ (Lemma 2). Thus $(t_i)$ and $(p_i)$ are
growth sequences. \hfill $\square$

\vspace{4mm}\noindent {\bf Theorem 3.}    {\sl Let $L$ be a cft
Lie algebra of finite depth and satisfying the growth conditions
(A.1) and (A.2). Set $\alpha =$ log index$\, L$. Then for some
$d$, $$\sum_{i=1}^{d-1}   \mbox{dim}\, L_{k+i}  = e^{(\alpha
+\varepsilon_k)k}\,, \hspace{1cm} \mbox{where}\, \varepsilon_k\to
0\, \mbox{as}\, k\to \infty \,.$$
 In
particular, this sum grows exponentially in $k$.}

\vspace{3mm}\noindent {\sl Proof.}  According to \cite{EHA} there
is a finitely generated sub Lie algebra $E \subset L$ such that
$\mbox{Ext}_{UL}^r(\bk , UL) \to \mbox{Ext}_{UE}^r(\bk,UL)$ is
non-zero.

\vspace{3mm}\noindent {\bf Lemma 5.}  {\sl The centralizer, $Z$,
of $E$ in $L$ is finite dimensional.}

\vspace{3mm}\noindent {\sl Proof.}  Since $E$ has finite depth, $Z
\cap E$ is finite dimensional \cite{rad}. Choose $k$ so $Z \cap E$
is concentrated in degrees $<k$. Suppose $x\in Z_{\geq k}$ has
even degree, and put $F = \bk x \oplus E$. Then
$\mbox{Ext}_{UF}(\bk, UF)\to \mbox{Ext}_{UE}(\bk, UF)$ is zero,
contradicting the hypothesis that the composite
$$\mbox{Ext}_{UL}^r(\bk, UL) \to \mbox{Ext}_{UF}^r(\bk, UL) \to
\mbox{Ext}_{UE}^r(\bk, UL)$$ is non-zero.

It follows that $Z_{\geq k}$ is concentrated in odd degrees, hence
an abelian ideal  in  $Z+E$. Again
$$\mbox{Ext}_{UL}^r(\bk, UL)\to \mbox{Ext}_{U(Z+E)}^r(\bk, UL) \to
\mbox{Ext}_{UE}^r (\bk, UL)$$ is non-zero. Thus
$\mbox{Ext}_{U(Z+E)}^r(\bk, UL)\neq 0$, $Z+E$ has finite depth and
every abelian ideal in  $Z+E$ is finite dimensional. \hfill $\square$

\vspace{5mm} Choose $d$ so that $E$ is generated in degrees $\leq
d-1$. As in the Proposition, set $\mbox{log index}\, L = \mbox{log
  index}\, L/Q = \alpha$.  If the theorem fails we can find a
sequence $s_i \to \infty$ such that $$\sum_{j=1}^{d-1}
\mbox{dim}\, L_{s_i+j} \leq e^{(\alpha - \beta)s_i}\,,
\eqno{(7)}$$ some $\beta >0$. Apply the Proposition to find growth
sequences $t_i$ and $p_i$ such that $t_i \leq s_i < p_i$ and
$p_i/t_i \to 1$.

We now use (7) to prove that $$\mbox{dim}\, (UE)_{(s_i-t_i,
s_i-t_i+d)} \, \mbox{dim}\, L_{(s_i,s_i+d)} < \mbox{dim}\,
L_{t_i}\,, \hspace{5mm} \mbox{$i$ large}\,. \eqno{(8)}$$ In fact
since $\limsup\, (\mbox{dim}\, (UL)_i)^{1/i} = e^\alpha$ it
follows that for some $\gamma>0$, $$\sum_{j=1}^{d-1}\mbox{dim}\,
(UE)_{j+k} \leq e^{\gamma (k+1)}\,, \hspace{1cm}\mbox{all
$k$}\,.$$ Thus it is sufficient to show that $$\gamma (s_i -
t_i+1) + (\alpha - \beta)s_i < (\alpha + \varepsilon_{t_i})t_i\,,
\hspace{1cm} \mbox{large $i$}\,,$$ where $\varepsilon_{t_i} \to 0$
as $i\to \infty$.  Write $s_i = \mu_it_i$; then $\mu_i \to 1$ and
the inequality reduces to the obvious
 $$\gamma/t_i + \gamma (\mu_i - 1) + (\alpha - \beta)\mu_i < (\alpha
+ \varepsilon_{t_i})\,, \hspace{1cm} \mbox{large}\, i\,.$$
Thus (8) is established.

Choose $s=s_i$, $t=t_i$ so that (8) holds and so that $Z_j = 0$,
$j\geq t$. Write $s-t = k$. The adjoint action of $UL$ in $L$
restricts to a linear map $$\left[ \displaystyle\oplus_{j=1}^{d-1}
(UE)_{k+j}\right] \otimes L_t \to
\displaystyle\oplus_{j=1}^{d-1}L_{s+j}$$ and it follows from (8)
that for some non-zero $x\in L_t$, $$(\mbox{ad}\, a)x = 0\,,
\hspace{1cm} a\in UE_{(k,k+d)}\,.$$

On the other hand, since $E$ is generated in degrees $\leq d-1$,
$(UE)_{> k} = UE \cdot (UE)_{(k,k+d)}$. Thus $(UE)_{>
k}\cdot x = 0$ and so $(UE)\cdot x $ is finite dimensional. A
non-zero element $y$ of maximal degree in $UE\cdot x$ satisfies
$$[a,y] = (\mbox{ad}\, a)(y) = 0\,, \hspace{1cm} a\in E\,,$$ i.e.
$y\in Z$ in contradiction to $Z_{\geq t} = 0$. This completes the
proof of the Theorem.\hfill $\square$

\section{Proof of Theorems 1 and 2}

\noindent{\sl Proof of Theorem 2:}   We show that $L_X$ satisfies
the hypothesis of Theorem 3. Since depth $L_X< \infty$
(\cite{rad}) and (A.2) holds by hypothesis we have only to
construct a quasi-geometric growth sequence $(r_i)$.

Let $\alpha = \mbox{log index}\, L_X$. Then $\alpha >0$ by
\cite{IHES}. Choose a sequence
 $$u_1 < u_2 < \cdots $$ such that
$(\mbox{dim}\, (L_X)_{u_i})^{1/u_i} \to e^\alpha$.

Next, suppose $\mbox{cat}\, X = m$ and put $a\left(\displaystyle\frac{1}{2(m+1)}\right)^{m+1}$. By starting the
sequence at some $u_j$ we may assume $\mbox{dim}\, (L_X)_{u_i} >
\frac{1}{a}$, all $i$. Thus the formula  in (\cite{IHES}, top of
page 189) gives a sequence $$u_i = v_0 < v_1 \cdots < v_ku_{i+1}$$ such that $v_{i+1} \leq 2(m+1)v_i$ and
$$\left(\mbox{dim}\,( L_X)_{v_j}\right)^{\frac{1}{v_j+1}}\geq \left[a\,
 \mbox{dim}\, (L_X)_{v_0}\right]^{\frac{1}{v_0+1}}\,, \hspace{5mm} j<k\,.$$
Since $v_0 = u_i$ and $ u_i\to \infty$ it follows that
$a\frac{1}{v_0+1}\to 1$ as $i\to \infty$.  Hence
 interpolating the sequences $u_i$ with the
sequences $v_j$ gives a quasi-geometric growth sequence $(r_j)$.
\hfill $\square$

\vspace{4mm}\noindent {\sl Proof of Theorem 1:}  A theorem of
Adams-Hilton \cite{AH} shows that $$UL_X = H_*(\Omega X;\mathbb Q)
= H(TV,d)$$ where $TV$ is the tensor algebra on $V$ and $V_i \cong
H_{i+1}(X;\mathbb Q)$. Thus $V$ is finite dimensional. Since $TV$
has a strictly positive radius of convergence so do $H(TV,d)$ and
$L_X$ : $$\mbox{log index}\, L_X <\infty\,.$$

Thus $X$ satisfies the hypotheses of Theorem 2. The fact that $d$
can be replaced by $n$  is proved by Lambrechts in \cite{La}.
\hfill $\square$

 \vspace{1cm}  Universit\'e Catholique de Louvain,
1348, Louvain-La-Neuve, Belgium

  University of Maryland, College Park, MD
20742-3281,USA

 Universit\'e d'Angers, 49045 Bd Lavoisier, Angers,
France

\end{document}